\documentclass[11pt]{article}
\usepackage{psfig,amssymb,latexsym}

\newtheorem{theorem}{Theorem}[section]

\newtheorem{prop}{Proposition}[section]

\def\slfrac#1#2{\hbox{\kern.1em %
 \raise.5ex\hbox{\the\scriptfont0 #1}\kern-.11em %
 /\kern-.15em\lower.25ex\hbox{\the\scriptfont0 #2}}}

\newcommand{\eqn}[1]{(\ref{#1})}
\newcommand{\hsp}{\hspace*{\parindent}}

\newcommand{\eeq}{\end{equation}}
\newcommand{\beql}[1]{\begin{equation}\label{#1}}
\newcommand{\bsq}{{\vrule height .9ex width .8ex depth -.1ex }}

\newcommand{\La}{\Lambda}
\newcommand{\la}{\lambda}

\newcommand{\Om}{\Omega}

\newcommand{\ZZ}{{\mathbb Z}}
\newcommand{\RR}{{\mathbb R}}

\newcommand{\sB}{{\cal B}}
\newcommand{\sC}{{\cal C}}

\newcommand{\sT}{{\cal T}}
\newcommand{\sA}{{\cal A}}

\makeatletter
\def\@sect#1#2#3#4#5#6[#7]#8{\ifnum #2>\c@secnumdepth
     \def\@svsec{}\else
     \refstepcounter{#1}\edef\@svsec{\csname the#1\endcsname.\hskip .75em }\fi
     \@tempskipa #5\relax
      \ifdim \@tempskipa>\z@
        \begingroup #6\relax
          \@hangfrom{\hskip #3\relax\@svsec}{\interlinepenalty \@M #8\par}%
        \endgroup
       \csname #1mark\endcsname{#7}\addcontentsline
         {toc}{#1}{\ifnum #2>\c@secnumdepth \else
                      \protect\numberline{\csname the#1\endcsname}\fi
                    #7}\else
        \def\@svsechd{#6\hskip #3\@svsec #8\csname #1mark\endcsname
                      {#7}\addcontentsline
                           {toc}{#1}{\ifnum #2>\c@secnumdepth \else
                             \protect\numberline{\csname the#1\endcsname}\fi
                       #7}}\fi
     \@xsect{#5}}
\def\@begintheorem#1#2{\it \trivlist \item[\hskip \labelsep{\bf #1\ #2.}]}

\def\plain{plain}\ifx\fmtname\plain\csname fi\endcsname
     
     \input docstrip
     \preamble

     Do not distribute the stripped version of this file.
     The checksum in the header refers to the documented version.

     \endpreamble
     \generateFile{here.sty}{t}{\from{here.doc}{}}
     \endinput
\fi
\ifcat a\noexpand @\let\next\relax\else\def\next{%
    \documentstyle[here,doc]{article}\MakePercentIgnore}\fi\next
\ifx\@Hxfloat\@Hundef\else\expandafter\endinput\fi
\let\@Hxfloat\@xfloat
\def\@xfloat#1[{\@ifnextchar{H}{\@HHfloat{#1}[}{\@Hxfloat{#1}[}}
\def\@HHfloat#1[H]{%
\expandafter\let\csname end#1\endcsname\end@Hfloat
\vskip\intextsep\vbox\bgroup\def\@captype{#1}\parindent\z@
\ignorespaces}
\def\end@Hfloat{\egroup\vskip \intextsep}

\makeatother

\catcode`\@=11
\renewcommand{\section}{
        \setcounter{equation}{0}
        \@startsection {section}{1}{\z@}{-3.5ex plus -1ex minus
        -.2ex}{2.3ex plus .2ex}{\large\bf}
        }
\catcode`@=12

\begin{document}
\begin{center}
{\Large {\bf  Universal Spectra and Tijdeman's Conjecture on Factorization 
of Cyclic Groups}} \\
\vspace{1.5\baselineskip}
{\em Jeffrey C. Lagarias} \\
\vspace*{0.5\baselineskip}
{\em S\'{a}ndor Szab\'{o}} \\
\vspace*{1.5\baselineskip}
(July 28, 2000) \\
\end{center}
\noindent{\bf ABSTRACT.} A spectral set $\Omega$ in $\RR^n$ is a set of finite Lebesgue measure 
such that
$L^2 ( \Omega )$ has an orthogonal basis of exponentials 
$\{ e^{2 \pi i \langle \la, x \rangle} : \la \in \La \}$ 
restricted to $\Omega$. Any such set $\La$
is called a spectrum for $\Omega$.
It is conjectured that every spectral set $\Omega$ tiles $\RR^n$ by 
translations.
A tiling set $\sT$ of translations has a {\em universal spectrum} $\La$ if 
every set
$\Omega$ that tiles $\RR^n$ by $\sT$ is a spectral set with spectrum $\La$.
Recently Lagarias and Wang showed that many periodic tiling sets $\sT$
have universal spectra. Their proofs used properties of factorizations of
abelian groups, and were valid for all groups for which a strong form
of a conjecture of Tijdeman is valid. However Tijdeman's original 
conjecture is not true in general, as follows from a construction 
of Szab\'o~\cite{Sz85}, and here we give a counterexample to 
Tijdeman's conjecture for the cyclic group of order 900. 
This paper formulates a new sufficient condition 
for a periodic tiling set to have a universal spectrum,
and applies it to show that the tiling sets in the given
counterexample do possess universal spectra. \\



\section{Introduction}
\hsp
A {\em spectral set} $\Omega$ in $\RR^n$ is a closed set of
finite Lebesgue measure such that $L^2 ( \Omega )$ has an orthogonal
basis of exponentials $\{ e^{2 \pi i \langle \la, x \rangle} : \la \in \La \}$
restricted to $\Omega$.
Any discrete set $\La$ in $\RR^n$ 
with this property is called a {\em spectrum} for $\Omega$
and $(\Omega, \La )$ is called a {\em spectral pair}.

In 1974 B. Fuglede \cite{Fu74} studied the problem of finding self-adjoint
commuting extensions of the operators
$i \frac{\partial}{\partial x_1} , \ldots, i \frac{\partial}{\partial x_n}$ 
inside $L^2 ( \Omega )$, and related the existence of such an
extension to $\Omega$ being a spectral set.
He formulated the following conjecture.

\paragraph{Spectral Set Conjecture.}
{\em Let $\Omega$ be a measurable set of $\RR^n$ with finite Lebesgue measure.
Then $\Omega$ is a spectral set if and only if $\Omega$ tiles $\RR^n$ 
by translations.
} \\

\noindent
Here $\Omega$ tiles $\RR^n$ with a set $\sT$ of translations
if $\Omega + \sT = \RR^n$ and, for $t, t' \in \sT$,
$$meas (( \Omega + t) \cap ( \Omega + t' )) =0 \quad\mbox{if}\quad
t \neq t' ~.
$$
Despite extensive study, this conjecture remains open in all dimensions,
in either direction.
For work on this problem see \cite{Fu74}, \cite{IP98},
\cite{J82}, \cite {JP92}, \cite{JP99},
\cite{K99},
\cite{LRW98}, 
\cite{LW97}, \cite{P96}.

Lagarias and Wang \cite{LW97} approached the spectral set conjecture in 
terms of tiling sets.
They proved that a large class of tiling sets $\sT$ had a 
{\em universal spectrum} $\La$ in the sense that every set 
$\Omega$ that tiles $\RR^n$ with the tiling set $\sT$
was a spectral set with the spectrum $\La$.
They considered periodic tiling sets of the form\footnote
{General periodic tilings can always be reduced to this form by a linear 
transformation of $\RR^n$.} 
$\sT = N_1\ZZ \times \dots \times N_n\ZZ  +  \sA$ with 
$\sA \subseteq \ZZ^n$, for some $N_1, \dots, N_n \in \ZZ$.
Their results were obtained by reduction to questions about factorizations~
\footnote{Factorizations of groups are 
defined in \S2.}
of abelian groups.
(A  connection to such factorizations was 
originally noted in Fuglede \cite{Fu74}.)
In the one-dimensional case, they formulated the following conjecture.

\paragraph{Universal Spectrum Conjecture.}
{\em Let $\sT = \ZZ + \frac{1}{m} A$ where $A \subseteq \ZZ$ 
reduced (mod~$m$) 
admits some factorization
$A \oplus B  = \ZZ / m \ZZ$.
Then $\sT$ has a universal spectrum of the form $m \ZZ + \Gamma$, 
where $\Gamma \subseteq \ZZ$.
} \\

In support of this conjecture,
they proved \cite[Theorem 1.2]{LW97} that if the cyclic group
$\ZZ_m$ has a property called the strong Tijdeman property, 
then any tile set
$\sT = \ZZ + \frac{1}{m} \sA$ has a universal spectrum of the form 
$\La = m \ZZ + \Gamma$ for some $\Gamma \subseteq \ZZ$.
They showed that the strong Tijdeman property holds for many cyclic
groups, and formulated the
conjecture that all cyclic groups have the strong Tijdeman property.
We do not define the strong Tijdeman property here (see \cite{LW97})
but observe that a necessary requirement for its truth for 
$G = \ZZ /m \ZZ$    
is the truth of the Tijdeman conjecture ($\bmod~m$) stated in \S2.
In addition to these results, Petersen and Wang \cite[Theorem 4.5]{PW99} 
also found other tiling  sets $\sT$ which have universal spectra.

Recently Coven and Meyerowitz~\cite{CM99} observed that a conjecture
equivalent to Tijdeman's Conjecture $(\bmod~m)$ had been made 
earlier, by Sands~\cite{Sa77} in 1977. Sands proved this conjecture
holds for all $m$ divisible by at most two distinct primes. 
Sands' conjecture $(mod ~)$ was disproved by Szab\'o~\cite{Sz85} in 1985,
by a direct construction which applies to certain integers $m$ divisible 
by three or more distinct primes. 
The smallest counterexample covered by that construction is 
$m = 2^3 3^3 5^2  = 5400.$ In \S2 we present a counterexample to
this conjecture for $m = 2^2 3^2 5^2 = 900,$ found using similar ideas to
the construction in Szab\'o~\cite{Sz85}.
The counterexamples disprove the  
Tijdeman Conjecture ($\bmod$ $m$) 
for such $m$, and this shows that the method for proving the 
existence of universal spectra 
given in \cite{LW97} does not work in general. 
In consequence, new  methods are needed 
to resolve the universal spectrum conjecture. 

In \S3 we present
a new sufficient condition for the existence of a universal spectrum
for a periodic tiling set in $\RR^n$ (Theorem~\ref{th31}). 
This criterion is easier to check computationally than a
necessary and sufficient condition given in \cite{LW97}.
It seems conceivable that the condition of Theorem~\ref{th31} is actually 
necessary and sufficient; this remains an unresolved question. 
If so it might prove useful in resolving the spectral set conjecture 
for periodic tiling sets.

In \S4 we apply the sufficient condition  of \S3 to
show that the tiling sets $\sT_A = \ZZ + \frac{1}{900} A$ and
$\sT_B = \ZZ + \frac{1}{900} B$ associated to the counterexample in \S2
do have universal spectra.
These sets give examples supporting the universal spectrum conjecture
which are not covered by the methods of \cite{LW97} and \cite{PW99}. 

\section{Counterexample to Tijdeman's Conjecture}
\hsp
A {\em factorization} $(A,B)$ of the finite cyclic group $G = \ZZ / m \ZZ$,
written
\beql{B201}
\bar{A} \oplus \bar{B} = G~,
\eeq
is one in which every element $g \in G$ has a unique representation
$$
g = \bar{a} + \bar{b} , \quad
\bar{a} \in \bar{A} \quad \mbox{and}\quad
\bar{b} \in \bar{B} ~.
$$
We write $\bar{A} = \{ \bar{a}_0, \ldots, \bar{a}_{n-1} \}$ and
$\bar{B} = \{ \bar{b}_0, \bar{b}_1, \ldots, \bar{b}_{l-1} \}$ with
$ln =m$.
Let $A = \{ a_0, \ldots, a_{n-1} \} \subseteq \ZZ$ be a lifting
of $\bar{A}$ to $\ZZ$, with $\bar{a}_i = a_i$ $(\bmod~m)$,
and similarly let $B$ be a lifting of $\bar{B}$.
Any factorization $(A,B)$ of $\ZZ /m \ZZ$ yields the {\em direct
sum decomposition} of $\ZZ$, as
\beql{B202}
A \oplus \sT_B = \ZZ ~,
\eeq
in which
\beql{B203}
\sT_B = m \ZZ + B
\eeq
is a periodic set with period $m \ZZ$.
There has been extensive study of the structure of factorization
of finite cyclic groups (more generally finite abelian groups)
and of direct sum decomposition (1.2) of the integers $\ZZ$,
the history of which can be found in Tijdeman \cite{Ti95}, see also
Coven and Meyerowitz \cite{CM99}. The original conjecture of
Tijdeman \cite[p. 266]{Ti95}  is as follows.

\paragraph{Tijdeman's Conjecture.}
{\em If $A \oplus \sT = \ZZ$ with $0 \in A \cap \sT$ and $A$ is a finite set
with $n$ elements, with g.c.d. $\{ a: a \in A \} =1$, then there exists
some prime factor $p$ of $n$ such that all elements of $\sT$ are 
divisible by $p$.
}\\ 

Haj\'{o}s \cite{Ha50} and de~Bruijn \cite{dB50} showed for any direct
sum decomposition $A + \sT = \ZZ$ where $|A|$ is finite, the 
infinite set $\sT$ is periodic and necessarily has the form 
\eqn{B203} and thus corresponds to some factorization \eqn{B201} 
of a cyclic group $\ZZ/ m \ZZ$ in which $|A|$ divides $m$.
This permits Tijdeman's conjecture to be reformulated as a series of 
conjectures for each finite cyclic group $\ZZ / m \ZZ$, as follows.
\paragraph{Tijdeman's Conjecture ($\bmod~m$).}
{\em
If $\bar{A} \oplus \bar{B} = \ZZ /m \ZZ$ with 
$\bar{0} \in \bar{A} \cap \bar{B}$ and if $\bar{A}$ lifts to set 
$A \subseteq \ZZ$ with $0 \in A$ and g.c.d. $\{a: a \in A \} =1$, then any
lifting $B \subseteq \ZZ$ with $0 \in B$ has \\
g.c.d. $\{ b: b \in B \} \neq 1$.
} \\

Tijdeman showed that this conjecture  holds for any $m$ for which
$\ZZ / m \ZZ$ is a ``good'' group in the sense of Haj\'{o}s \cite{Ha50} and
de~Bruijn \cite{dB50}.
The complete list of cyclic ``good'' groups are known to be those 
of order $p^n$
$(n \ge 1)$, $pq$, $pqr$, $p^n q$ $(n > 1)$,
$p^2 q^2$, $p^2 qr$ and $pqrs$, where $p$, $q$, $r$ and $s$ are distinct
primes, c.f. \cite{Ti95} for references. As indicated in the introduction,
this conjecture was actually made earlier by Sands\cite{Sa77}, who proved
it also holds in the cases $m = p^n q^k$ $(n,k \geq 1)$. 

In 1985 Szab\'o~\cite{Sz85} gave a construction which gave counterexamples
to Sands' conjecture for certain $m$ divisible by three or more distinct 
primes. Szab\'o actually constructs
sets $A \subset \ZZ$ with $0 \in A$ and $\gcd(A) = 1$, which tile
the integers, whose members
are not uniformly distributed $\bmod ~k $ for any $ k \geq 2$.
Coven and Meyerowitz\cite[Lemma 2.5]{CM99} observe that it follows
that any tiling set $C \subset \ZZ$ for $A$ with $0 \in C$ must
have $\gcd (C) = 1,$ hence cannot be contained in any subgroup
of $\ZZ /m \ZZ$, where  $m$ is the minimal period of $C$.
The Szab\'o construction applies to $m = m_1 m_2 ... m_r$ with $r \geq 3$
in which each $m_i = u_i v_i$ with the $m_i$ pairwise relatively
prime, $u_i$ is the smallest prime dividing $m_i$ and each $v_i \geq 4.$
The smallest $m$ satisfying these conditions is $m = 2^3 3^3 5^2 = 5400.$
Here we present a counterexample for $m = 2^2 3^2 5^2 = 900$ which was
found using similar ideas. 

\begin{theorem}\label{th21}
Tijdeman's Conjecture $(\bmod~900 )$ is false.
\end{theorem}
\paragraph{Proof.}
We take sets $A$ and $B$ with $|A| = |B| =30$.
The set
\beql{B204}
A = \{ 0, 36, 72, 108, 144\} \oplus \{0, 100, 200 \} \oplus \{ 0, 225 \} ~.
\eeq
It has $g.c.d. ~ \{ a: a \in A \} = 1,$
since g.c.d.$\{2^2 3^2, 2^2 5^2 , 3^2 5^2 \} =1$.
We choose
\beql{B205}
B = \left\{
\begin{array}{rrrrrrrrrr}
0 & 30 & 60 & 126 & 180 & 210 & 220 & 240 & 300 & 306 \\
330 & 360 & 375 & 390 & 480 & 486 & 510 & 520 & 540 & 570 \\
660 & 666 & 690 & 750 & 780 & 820 & 825 & 840 & 846 & 870
\end{array}
\right\} ~.
\eeq
We claim that $\bar{A} \oplus \bar{B} = \ZZ / 900 \ZZ$ is a direct sum.
This can be verified by a calculation\footnote{More generally one
can consider the construction in Szab\'o\cite{Sz85}.} 
(a short computer program).
However
$$g.c.d. ~ \{ b : b \in B_0\}  =1,$$
by considering the values 126, 220 and 375. 
Thus Tijdeman's conjecture ($\bmod ~900$) is false.~~~$\bsq$

\paragraph{Remark.} Haj\'{o}s \cite{Ha50} advanced a weaker 
conjecture concerning direct sum 
decompositions of a cyclic group $G$, 
which is that every factorization $\bar{A} \oplus \bar{B} = G$ 
is quasiperiodic. We say that
a factorization is {\em quasiperiodic} if one of $\bar{A}$ or $\bar{B}$,
say $\bar{B}$, can be partitioned into disjoint subsets
$\{ \bar{B}_1, \ldots, \bar{B}_m \}$ such that there is a subgroup 
$H = \{h_1, \ldots, h_m \}$ of $G$ with
$$\bar{A} + \bar{B}_i = \bar{A} + \bar{B}_1 + h_i , \quad
1 \le i \le m ~.
$$
The example $(\bar{A}, \bar{B})$ above for $G = \ZZ / 900 \ZZ$ is 
quasiperiodic.
The choices are $H = \{ 0, 300, 600 \}$ and
\begin{eqnarray*}
\bar{B}_1 & = & \{
\begin{array}{rrrrrrrrrr}
\,~0 & 126 & 180 & 306 & 360 & 486 & 540 & 666 & 820 & 846
\end{array}
\} \\
\bar{B}_2 & = & \{
\begin{array}{rrrrrrrrrr}
30 & 210 & 220 & 300 & 390 & 480 & 570 & 660 & 750 & 840
\end{array}
\} \\
\bar{B}_3 & = & \{
\begin{array}{rrrrrrrrrr}
60 & 240 & 330 & 375 & 510 & 520 & 690 & 780 & 825 & 870
\end{array}
\} ~.
\end{eqnarray*}
The quasiperiodicity conjecture remains open.

\section{Criterion for Universal Spectrum}
\hsp
We formulate a sufficient condition for a universal spectrum for a periodic
tiling set in $\RR^n$, which is simpler to check than the necessary and
sufficient condition given in  \cite[Theorem 1.1]{LW97}. 

Given a finite set $\sB \subseteq \RR^n$ let $f_\sB (\la )$ denote the 
exponential polynomial
\beql{301}
f_\sB ( \la ) := \sum_{b \in \sB} e^{2 \pi i <\la, b >}
\eeq
and let $Z( f_\sB )$ denote its set of real zeros, i.e.
\beql{302}
Z( f_\sB ) = \{ \la \in \RR^n : f_\sB (\la) =0 \} ~.
\eeq
We recall the following criterion for a set $\Lambda$ to be
a spectrum, taken from \cite{LW97}.
\begin{prop}\label{pr31}
Let $\Omega = [ 0, \frac{1}{N_1} ] \times \dots \times [ 0, \frac{1}{N_n}]
 + \sB$ where 
$\sB \subseteq \frac{1}{N_1} \ZZ \times \dots \times 
\frac{1}{N_n} \ZZ$ is a finite set.
Suppose that $\Gamma \subseteq \ZZ^n$ is a set of distinct residue classes
$( \bmod ~N_1 \ZZ \times \dots \times N_n \ZZ)$, i.e.
$( \Gamma - \Gamma ) \cap (N_1 \ZZ \times \dots \times N_n \ZZ) = \{ {\bf 0} \}$.
Then $\Lambda = (N_1 \ZZ \times \dots \times N_n\ZZ) + \Gamma$ 
is a spectrum for $\Omega$ if and only if
$| \Gamma | = | \sB |$ and
\beql{303}
\Gamma - \Gamma \subseteq Z (f_\sB ) \cup \{ {\bf 0}\} ~.
\eeq
\end{prop}

\paragraph{Proof.}
This is \cite[Theorem 2.3]{LW97}, after a linear
rescaling of Euclidean space $\RR^n$ of $\Omega$ by a factor 
$(\frac{1}{N_1}, ..., \frac{1}{N_n})$ and 
a corresponding dilation of Fourier space by a factor $(N_1,...,N_n)$.~~~$\bsq$

\begin{theorem}\label{th31}
$($Universal Spectrum Criterion$)$.
Let $\sT = \ZZ^n + \sA$ where 
$\sA \subseteq \frac{1}{N_1} \ZZ \times  \dots \times \frac{1}{N_n} \ZZ$
is a finite set, 
and suppose there exists some set $\Omega$ such that $\Omega$ tiles $\RR^n$
by translations using the tiling set $\sT$. Consider a set  
$\Lambda = (N_1 \ZZ \times \dots \times N_n\ZZ) + \Gamma$ with 
$\Gamma \subseteq \ZZ^n$ such that 
the residue classes 
$\Gamma~(mod~ N_1 \ZZ \times \dots \times N_n \ZZ)$ are all distinct, i. e.
\beql{303a}
(\Gamma - \Gamma ) \cap (N_1 \ZZ \times \dots \times N_n \ZZ) = \{ 0 \} 
\eeq
Then $\Lambda$ is a universal spectrum for  $\sT$ provided that 
$| \Gamma | = \frac {N_1 N_2 ... N_n}{|\sA|}$,  and 
\beql{304}
(\Gamma - \Gamma ) \cap Z ( f_\sA ) = \emptyset ~.
\eeq
\end{theorem}
\paragraph{Proof.}
By Theorem 3.1 of  \cite{LW97} a set $\Omega$ tiles $\RR^n$ by translations
with a periodic tiling set $\sT$ if and only if there exists some finite set
$\sB \subseteq \frac{1}{N_1} \ZZ \times  \dots \times \frac{1}{N_n} \ZZ$
 giving a factorization
$$\sA \oplus \sB = \ZZ_{N_1} \times \dots \times \ZZ_{N_n} = (\frac{1}{N_1} \ZZ / \ZZ)
\times \dots \times (\frac{1}{N_n}\ZZ / \ZZ),$$
in which case $\tilde{\Omega} = [0, 1]^n + \sB$ also has $\sT$ as a tiling set.

By Theorem 1.1 of \cite{LW97} it suffices to verify that $\Lambda$ 
is a spectrum for each set
$$\Omega_{\sB} = [ 0, \frac{1}{N_1}]\times \dots \times [ 0, \frac{1}{N_n} ] + \sB$$
where $$\sA \oplus \sB = (\frac{1}{N_1}\ZZ/ \ZZ) \times \dots \times 
(\frac{1}{N_n} \ZZ  / \ZZ)$$ is a direct sum decomposition.
This property shows that
$$\Omega_{\sB} + \sA = [ 0, \frac{1}{N_1}]\times \dots 
\times [ 0, \frac{1}{N_n} ] + \sB + \sA$$ 
is a fundamental domain for the 
$n-$torus $\RR^n / \ZZ^n.$
It follows that the Fourier transform
$$\hat{\chi}_{\Omega + \sA} ( \la ) = \int_{\Omega + \sA} e^{2 \pi i <\la, x>} dx$$
for $\la = k \in \ZZ^n$ satisfies
\beql{305}
\hat{\chi}_{\Omega + \sA} (k) = \left\{
\begin{array}{lll}
1 & \mbox{for} & k={\bf 0} , \\
0 & \mbox{for} & k \in \ZZ^n \setminus \{0\} ~.
\end{array}
\right.
\eeq
Now
\begin{eqnarray*}
\hat{\chi}_\Omega ( \la ) & = & \int_\Om e^{2 \pi i <\la, x>} dx =
\sum_{b \in \sB} e^{2 \pi i <\la,  b>} \prod_{j=1}^n
\int_0^{\frac{1}{N_j}} e^{2 \pi i \la_j x_j} dx_j \\
& = & f_\sB (\la ) \prod_{j=1}^n \frac{e^\frac{2 \pi i \la_j}{N_j} -1}{2 \pi i \la_j} ~,
\end{eqnarray*}
which gives
\beql{306}
\hat{\chi}_{\Om + \sA} ( \la ) = f_\sA ( \la ) f_\sB (\la)
\prod_{j=1}^n \left( \frac{\sin  \frac{\pi\la_j}{N_j}}{\pi \la_j} \right)
e^{\pi i \frac{\la_j}{N_j}} ~.
\eeq
Comparing \eqn{305} and \eqn{306} yields
\beql{307}
f_\sA (k) f_\sB (k) =0 \quad\mbox{if} \quad k \in \ZZ^n \setminus 
(N_1\ZZ \times \dots \times N_n\ZZ) ~.
\eeq
Thus we obtain
\beql{308}
f_\sB (k) = 0 \quad \mbox{if} \quad k \not\in Z(f_\sA ) \cap 
(\ZZ^n \setminus (N_1\ZZ \times \dots \times N_n\ZZ)) ~.
\eeq
By hypothesis $\La = (N_1\ZZ \times \dots \times N_n\ZZ)  + \Gamma$ with 
$\Gamma \subseteq \ZZ^n$, $| \Gamma | = \frac {N_1 \dots N_n}{| \sA |}$,
\beql{309}
(\Gamma - \Gamma ) \cap (N_1\ZZ \times \dots \times N_n\ZZ)  = \{{\bf 0}\} ~,
\eeq
and $(\Gamma - \Gamma ) \cap Z( f_\sA ) = \emptyset$.
Thus $|\Gamma | = | \sB |$.
We claim that
\beql{310}
\La - \La \subseteq Z (f_\sB ) \cup \{{\bf 0}\} ~.
\eeq
This claim holds since $\La - \La \subseteq \ZZ^n$, 
then noting that $Z( f_\sB )$ contains all points of
$\ZZ^n \setminus (N_1\ZZ \times \dots \times N_n\ZZ)$ not in $Z(f_\sA )$ 
by \eqn{308}, while \eqn{309} takes care of
points in $N_1\ZZ \times \dots \times N_n\ZZ$.
Now Proposition \ref{pr31} shows that $\Omega$ has $\La$ as a spectrum, 
and the theorem follows.~~~$\bsq$

\paragraph{Remarks.}
(1). The main hypothesis in Theorem~\ref{th31} is \eqn{304},
which requires determining the finite
set $\Gamma - \Gamma \subseteq \ZZ^n$ and evaluating $f_\sA$ at 
these points. This
condition is computationally simpler to check than the criterion of 
\cite[Theorem 1.1]{LW97}, which requires
determining all the complementing sets $\sB$ to $\sA$.

(2). It seems conceivable that the sufficient conditon
 of Theorem \ref{th31} might  
also be a necessary condition for a universal spectrum of the given 
form $\La = m \ZZ + \Gamma$ with
$\Gamma \subseteq \ZZ$.
To show this one would have to show that for each integer vector
${\bf k} \in Z (f_\sA ) \cap ( \ZZ^n \setminus N_1\ZZ \times \dots \times N_n\ZZ )$ 
there exists some set $\sB$
with $\sA \oplus \sB = (\frac{1}{N_1} \ZZ / \ZZ) \times \dots \times 
(\frac{1}{N_n} \ZZ / \ZZ) $ such that
${\bf k} \not\in Z(f_\sB )$.

\section{Universal Spectra}

We apply the criterion of Theorem~\ref{th31} 
to show that the the tiling sets $\ZZ + \frac{1}{900} A$ and 
$\ZZ + \frac{1}{900} B$ 
associated to this counterexample  in \S2 both have universal spectra.

\begin{theorem}\label{th32}
Let $\sT_A = \ZZ + \frac{1}{900} A$ with $A = \{0,36,72,108,144\} \oplus \{ 0,100,200\} \oplus \{0, 225\}$.
Then $\sT_A$ has the universal spectrum $\La_A = 900 \ZZ + A$.
\end{theorem}

\paragraph{Proof.}
We apply Theorem \ref{th31} with $n=1$ and $N_1=900$.
Then $\sA = \frac{1}{900} A$ and $\sB = \frac{1}{900} B$ with $B$ given in
\eqn{B205}
gives a factorization $\sA \oplus \sB = \frac{1}{900} \ZZ / \ZZ$.
A calculation gives
\beql{311}
f_\sA ( \la ) = \left(
\frac{1- e^{\frac{2 \pi i \la}{5}}}{1- e^{\frac{2 \pi i \la}{25}}}
\right)
\left(
\frac{1- e^{\frac{2 \pi i\la}{3}}}{1- e^{\frac{2 \pi i\la}{9}}} \right)
\left( \frac{1- e^{\frac{2 \pi i \la}{2}}}{1- e^{\frac{2 \pi i \la}{4}}}
\right) ~.
\eeq
It follows that the set $Z( f_\sA ) \subseteq \ZZ$ consists of all integers
$k$ such that one or more of the following three conditions hold:
\begin{itemize}
\item[(i)]
5 divides $k$ and 25 doesn't divide $k$.
\item[(ii)]
3 divides $k$ and 9 doesn't divide $k$.
\item[(iii)]
2 divides $k$ and 4 doesn't divide $k$.
\end{itemize}
The set $\Lambda_A = 900 \ZZ + A$ has $\Lambda \subseteq \ZZ$ and
$(A-A) \cap 900 \ZZ = \{0\}$, since all $a \in A$ have $0 \le a \le 900$
and are distinct.
Also $|A| = 130 = \frac{900}{| \sA |}$.
To apply Theorem \ref{th31} it remains to verify
\beql{312}
(A-A) \cap Z( f_\sA ) = \emptyset ~.
\eeq
While $a_i \in A$ for $i=1,2$ as
$$a_i = 36 k_i + 100 l_i + 225 m_i$$
with $0 \le k_i \le 4$, $0 \le l_i \le 2$ and $0 \le m_i \le 1$.
Then
$$a_1 - a_2 = 36 (k_1- k_2) + 100 (l_1 - l_2 ) + 225 (m_1 - m_2 )$$
with $|k_1 - k_2 | < 5$, $|l_1 - l_2 | < 3$ and $|m_1 - m_2 | < 2$.
Thus if 5 divides $a_1 -a_2$, then 5 divides $k_1-k_2$ so $k_1 = k_2$,
and we conclude 25 divides $a_1 - a_2$.
By similar arguments if 3 divides $a_1 - a_2$ then 9 divides $a_1 - a_2$, while if 2 divides $a_1 - a_2$ then 4 divides
$a_1 - a_2$.
Thus none of (i), (ii), (iii) hold, and \eqn{312} follows.~~~$\bsq$

\paragraph{Remark.} The proof of Theorem \ref{th32} 
easily generalizes to 
the sets $A$ as appearing in the general construction of Szab\'o~\cite{Sz85}:
All such tiling sets $\ZZ + \frac{1}{m}A$ 
have a universal spectrum.

\begin{theorem}\label{th33}
Let $\sT_B = \ZZ + \frac{1}{900} B$ with $B$ given by {\rm \eqn{B205}}.
Then $\sT_B$ has the universal spectrum $\La_B = 900 \ZZ + B$.
\end{theorem}

\paragraph{Proof.}
We do not have a conceptual proof of this result;
however the conditions of Theorem \ref{th31} can be verified by a direct
calculation (on the computer.)
A key fact is that $Z(f_\sB ) \cap \ZZ \subseteq \ZZ \setminus 900 \ZZ$
and that the
complement of the set $\ZZ( f_\sB ) \cap \ZZ$ in $\ZZ$ is exactly\footnote
{To verify this, by {\rm \eqn{309}} it suffices to check that $f_\sB (k) \neq 0$ whenever $k \in Z( f_\sA )$.}
$Z(f_\sA ) \cup 900 \ZZ$.
Because this fact holds, one can prove $\La_B$ is a universal spectrum
for $\sT_B$ by checking that
$$B - B \subseteq Z (f_\sA ) \cap \{0\} ~.$$
Since $Z(f_\sA )$ is given by conditions (i)--(iii) in the proof of
Theorem \ref{th32}, it suffices to verify that every nonzero element
of $B-B$ satisfies one of (i)--(iii).
This can be done by hand.~~~$\bsq$

\paragraph{Remarks.}
In Theorem \ref{th32} and Theorem \ref{th33} the 
universal spectrum
$\Lambda$ exhibited is a scaled version of the tiling set $\sT$.
This fact is special to these examples, and cannot hold in
general. 
Any tiling set $\sT = \ZZ + \frac{1}{m} \sC$ with the property that
$\Lambda = m \ZZ +\sC$ is a universal spectrum must have
$m = |\sC|^2$, by Proposition \ref{pr31}.

(2). The proof of Theorem \ref{th32} easily generalizes to apply to
the sets $A$ appearing in the general construction of Szab\'o~\cite{Sz85}:
All such tiling sets $\ZZ + \frac{1}{m}A$ 
have a universal spectrum. 

\paragraph{Acknowledgment.}
The first author is indebted to J. A. Reeds for helpful computations, and
especially to M. Szegedy for computations and references.

{\rm AT\&T Labs--Research, Florham Park, NJ 07932-0971, USA} \\
{\em E-mail address:}~{\tt jcl@research.att.com} \\

{\rm Dept. of Mathematics, Univ. of Bahrain, P.O. Box 30238 ,
Isa Town, BAHRAIN }\\
{\em E-mail address:}~{\tt sszabo7@hotmail.com} 

\clearpage
\begin{thebibliography}{99}
\bibitem{dB50}
N. G. de Bruijn,
On bases for the set of integers, 
{\em Publ. Math. Debrecen} {\bf 1} (1950), 232--242.

\bibitem{dB55}
N. G. de Bruijn,
On the factorization of cyclic groups, 
{\em Indag. Math.} {\bf 17} (1955), 370--377.

\bibitem{CM99}
E. Coven and A. Meyerowitz,
Tiling the integers with translates of one finite set,
J. Algebra {\bf 212} (1999), 161--174.

\bibitem{Fu74}
B. Fuglede,
Commuting self-adjoint partial differential operators and a group theoretic
problem,
{\em J. Funct. Anal.} {\bf 16} (1974), 101--121.

\bibitem{Ha50}
G. Haj\'{o}s,
Sur la factorisation des groupes ab\'{e}liens,
{\em Casopis P\v{e}st Mat. Fys.} {\bf 74} (1950), 157--162.

\bibitem{IP98}
A. Iosevich and S. Pedersen,
Spectral and Tiling Properties of the Unit Cube, 
Inter. Math.
Res. Notices {\bf 16} (1998), 819--828.

\bibitem{J82}
P. E. T. Jorgensen,
Spectral theory for finite volume domains in $\RR^n$,
{\em Adv. Math.} {\bf 44} (1982), 105--120.

\bibitem{JP92}
P. E. T. Jorgensen and S. Pedersen,
Spectral theory for Borel sets in $\RR^n$ of finite measure,
{\em J. Funct. Anal.} {\bf 107} (1992), 72--104.

\bibitem{JP99}
P. E. T. Jorgensen and S. Pedersen,
Spectral Pairs in Cartesian Coordinates,
{\em J. Fourier. Anal. Appl.} {\bf 5} (1999). 285--302.

\bibitem{K99}
M. Kolountzakis, 
Packing, tiling, orthogonality and completeness,
Bull. London Math. Soc., to appear; arXiv math.CA/9904066.

\bibitem{LRW98}
J. C. Lagarias, J. A. Reeds and Y. Wang,
Orthonormal bases of exponentials for the $n$-cube,
Duke Math. J. {\bf 103} (2000), 25--37.

\bibitem{LW96}
J. C. Lagarias and Y. Wang,
Tiling the line with translates of one tile,
{\em Invent. Math.} {\bf 124} (1996), 341--365.

\bibitem{LW97}
J. C. Lagarias and Y. Wang,
Spectral Sets and Factorizations of Finite Abelian Groups,
{\em J. Funct. Anal.} {\bf 145} (1997), 73--98.

\bibitem{P96}
S. Pedersen,
Spectral sets whose spectrum is a lattice with a base,
{\em J. Funct. Anal.} {\bf 141} (1996), 496--509.

\bibitem{PW99}
S. Pedersen and Y. Wang,
Universal Spectra, Universal Tiling Sets and the Spectral Set Conjecture, 
preprint.

\bibitem{Sa77}
A. D. Sands,
On Keller's Conjecture for cetain cyclic groups,
Proc. Edinburgh Math. Soc. {\bf 22} (1977), 17--21.

\bibitem{Sz85}
S. Szab\'o,
A type of factorization of finite abelian groups,
Discrete Math. {\bf 54} (1985), 121--124.

\bibitem{Ti95}
R. Tijdeman,
Decomposition of the integers as a direct sum of two subsets,
in: ``Number Theory Seminar Paris 1992--1993'' (S. David, Ed.),
pp.~261--276.
Cambridge Univ. Press, Cambridge, 1995.
\end{thebibliography}
\end{document}